\documentclass[10pt, a4paper]{amsart}

\usepackage{amssymb}
\usepackage{amscd}
\usepackage{enumerate}
\usepackage{hyperref}
\usepackage[utf8]{inputenc}
\usepackage{newunicodechar}
\usepackage{varioref}
\usepackage[arrow,curve,matrix]{xy}

\usepackage{colortbl}
\usepackage{graphicx}

\usepackage{mathrsfs}

\usepackage{url}

\definecolor{linkred}{rgb}{0.7,0.2,0.2}
\definecolor{linkblue}{rgb}{0,0.2,0.6}

\numberwithin{figure}{section}

\usepackage[top=30truemm,bottom=30truemm,left=25truemm,right=25truemm]{geometry}

\newcommand\danger{
 \makebox[1.4em][c]{%
 \makebox[0pt][c]{\raisebox{.1em}{\tiny!}}%
 \makebox[0pt][c]{
 $\bigtriangleup$}}}

\newunicodechar{α}{\ensuremath{\alpha}}
\newunicodechar{β}{\ensuremath{\beta}}
\newunicodechar{χ}{\ensuremath{\chi}}
\newunicodechar{δ}{\ensuremath{\delta}}
\newunicodechar{ε}{\ensuremath{\varepsilon}}
\newunicodechar{Δ}{\ensuremath{\Delta}}
\newunicodechar{η}{\ensuremath{\eta}}
\newunicodechar{γ}{\ensuremath{\gamma}}
\newunicodechar{Γ}{\ensuremath{\Gamma}}
\newunicodechar{ι}{\ensuremath{\iota}}
\newunicodechar{κ}{\ensuremath{\kappa}}
\newunicodechar{λ}{\ensuremath{\lambda}}
\newunicodechar{Λ}{\ensuremath{\Lambda}}
\newunicodechar{μ}{\ensuremath{\mu}}
\newunicodechar{ν}{\ensuremath{\nu}}
\newunicodechar{ω}{\ensuremath{\omega}}
\newunicodechar{Ω}{\ensuremath{\Omega}}
\newunicodechar{π}{\ensuremath{\pi}}
\newunicodechar{φ}{\ensuremath{\phi}}
\newunicodechar{Φ}{\ensuremath{\Phi}}
\newunicodechar{ψ}{\ensuremath{\psi}}
\newunicodechar{Ψ}{\ensuremath{\Psi}}
\newunicodechar{ρ}{\ensuremath{\rho}}
\newunicodechar{σ}{\ensuremath{\sigma}}
\newunicodechar{Σ}{\ensuremath{\Sigma}}
\newunicodechar{τ}{\ensuremath{\tau}}
\newunicodechar{θ}{\ensuremath{\theta}}
\newunicodechar{Θ}{\ensuremath{\Theta}}
\newunicodechar{ξ}{\ensuremath{\xi}}
\newunicodechar{ζ}{\ensuremath{\zeta}}

\newcommand{\sH}{\mathcal{H}}
\newcommand{\sI}{\mathcal{I}}

\newcommand{\sK}{\mathcal{K}}

\newcommand{\sM}{\mathcal{M}}

\newcommand{\OO}{\mathcal{O}}

\newcommand{\sZ}{\mathcal{Z}}
\newcommand{\ol}[1]{\overline{#1}}

\newcommand{\GG}{\mathbb{G}}

\newcommand{\Zar}{\mathrm{Zar}}
\newcommand{\Nis}{\mathrm{Nis}}
\newcommand{\et}{\mathrm{et}}
\newcommand{\cdh}{\mathrm{cdh}}
\newcommand{\procdh}{\mathrm{procdh}}

\newcommand{\Kan}{\mathrm{Kan}}

\newcommand{\Gal}{\operatorfont{Gal}} 
\newcommand{\Hdg}{\operatorfont{Hdg}} 
\newcommand{\Crys}{\operatorfont{Crys}}  

\newcommand{\Sch}{\mathrm{Sch}}
\newcommand{\Et}{\operatorname{Et}}
\newcommand{\Sm}{\operatorname{Sm}}

\newcommand{\qcqs}{\mathrm{qcqs}}
\newcommand{\Spt}{\mathrm{Spt}}
\newcommand{\PSh}{\mathrm{PSh}}
\newcommand{\Shv}{\mathrm{Shv}}
\newcommand{\SmCor}{\mathrm{SmCor}}
\newcommand{\Cor}{\mathrm{Cor}}
\newcommand{\SmProj}{\mathrm{SmProj}}
\newcommand{\DM}{\mathrm{DM}}
\newcommand{\MS}{\mathrm{MS}}

\newcommand{\Spec}{\mathrm{Spec}}

\newcommand{\red}{\mathrm{red}}

\theoremstyle{theorem}
\newtheorem{theo}{Theorem}

\newtheorem{prop}[theo]{Proposition}

\newtheorem*{disc}{Disclaimer}

\theoremstyle{definition}
\newtheorem{defi}[theo]{Definition}

\newtheorem{obse}[theo]{Observation}
\newtheorem{rema}[theo]{Remark}

\newtheorem{exam}[theo]{Example}

\newtheorem{cons}[theo]{Construction}

\DeclareSymbolFontAlphabet{\scr}{rsfs}

\newcommand{\CC}{\mathbb{C}}
\newcommand{\FF}{\mathbb{F}}

\newcommand{\LL}{\mathbb{L}}
\newcommand{\HH}{\mathbb{H}}
\newcommand{\NN}{\mathbb{N}}
\newcommand{\PP}{\mathbb{P}}
\newcommand{\QQ}{\mathbb{Q}}

\newcommand{\ZZ}{\mathbb{Z}}
\renewcommand{\AA}{\mathbb{A}}

\newcommand{\m}{\mathfrak{m}}

\DeclareMathOperator{\id}{id}

\DeclareMathOperator{\dlog}{dlog}

\newcommand{\cM}{\mathcal{M}}

\newcommand{\eff}{\operatorname{eff}}
\newcommand{\op}{\operatorname{op}}

\newcommand{\dR}{\operatorname{dR}}
\newcommand{\crys}{\operatorname{crys}}
\newcommand{\EM}{\operatorname{EM}}
\newcommand{\Cone}{\operatorname{Cone}}

\title{Two new motivic complexes for non-smooth schemes}
\author{Shane Kelly}

\begin{document}

\maketitle

\begin{center}
Expanded notes from a talk at the RIMS Workshop, \\
Algebraic Number Theory and Related Topics \\
December 13th, 2023
\end{center}

\begin{disc}
These are expanded notes for a talk at a general algebraic number theory conference. The talk was very much aimed at non-experts. As such, for the most part the machinery is kept as minimal as possible, e.g., varieties over a field instead of quasi-compact quasi-separated schemes, triangulated categories instead of stable $\infty$-categories, spectral sequences instead of a filtered spectra, etc. Consequently, many distinctions, resp. hypotheses, that become important in the general setting are collapsed, resp. automatic. We made some attempt at addressing this in Remark~\ref{rema:ZEMrem} and Remark~\ref{rema:cdhZ}, but the expert reader is advised to consult the original articles \cite{EM23} and \cite{KS24}.
\end{disc}

For smooth schemes $X$ over a field, Voevodsky's motivic cohomology $H^i(X, \ZZ(j))$ (or equivalently,%
\footnote{Bloch's higher Chow groups are a Borel-Moore homology theory, rather than a cohomology theory. However, over smooth schemes there is a good duality theory, so up to shift and twist cohomology can be identified with Borel-Moore homology.} %
Bloch's higher Chow groups $CH^i(X, 2i{-}j)$) satisfies many of the properties predicted by Beilinson and Lichtenbaum, \cite{Bei82}, \cite{Bei87}, \cite{Lic83}. For example, it sits in an Atiyah-Hirzebruch spectral sequence
\begin{equation}
\HH_{\Zar}^{i-j}(X, \ZZ(-j)) \implies K_{-i-j}(X)
\end{equation}
calculating algebraic $K$-theory, and with finite coefficients prime to the characteristic, it coincides with étale cohomology in the range $i \leq j$.

For non-smooth schemes over a field admitting an appropriate resolution of singularities Voevodsky's cdh motivic cohomology\footnote{See Remark~\ref{rema:cdhZ} for clarification about what this means.} $\ZZ(n)^{\cdh}$ is still useful. This is demonstrated by Suslin and Voevodsky's work on the Bloch-Kato conjecture, \cite[\S 5]{SV00}. However, we lose some desirable features such as the Atiyah-Hirzebruch spectral sequence discussed in Example~\ref{exam:cdhBad}.

In this talk we discussed two approaches to correcting this that were developed concurrently. The one we discuss first is joint work of the author with Shuji Saito, \cite{KS24}. It is motivated by Voevodsky's topological 
approach to motivic cohomology. The idea is to coarsen the cdh topology slightly in order to preserve nilpotents. The coarsening we have chosen%
\footnote{Of course others options are available, see \cite{Kel24}.} %
is related via Proposition~\ref{prop:MVprocdh} to procdh excision studied by Morrow, \cite{Mor16}, and used in Kerz, Strunk, Tamme's proof of Weibel's conjecture, \cite{KST18}. We call it the \emph{procdh topology}. The resulting complexes are written 
\[ \ZZ(n)^{\procdh}. \]

The one we discuss second, developed independently and at the same time, is due to Elmanto and Morrow, \cite{EM23}. It is motivated by trace methods in algebraic $K$-theory. The idea is to repair the damage of cdh sheafifying by glueing back in the desired missing part. Explicitly, there is some appropriate map $\ZZ(n) \to \ZZ(n)^{TC}$ of complexes induced by the trace map $K(X) \to TC(X)$ and one can define 
\[ \ZZ(n)^{\EM} \]
as the $\Cone(\ZZ(n)^{\cdh} \oplus \ZZ(n)^{TC} \to \ZZ(n)^{TC, \cdh})[-1]$. The identification of $\ZZ(n)^{TC}$ and $\ZZ(n)^{TC, \cdh}$ depends on the characteristic; see Section~\ref{sec:4} for details.

These are qualitatively different approaches. The procdh approach is similar to Voevodsky's original constructions, e.g., the one that used the $h$-topology, \cite{Voe96}. It is automatically well-defined, even in mixed characteristic, and since it is defined via a universal property, there is automatically a map from $\ZZ(n)^{\procdh}$ to other potential candidates for motivic cohomology satisfying procdh excision. On the other hand the trace methods approach is more convenable to computations. It is difficult to concretely calculate $\ZZ(n)^{\procdh}$ without the comparison to $\ZZ(n)^{\EM}$.

Both of these approaches only give a cohomology, disembodied from a category of motives. In the last section we mention recent work of Annala, Iwasa and Annala, Iwasa, Hoyois one hopes will lead to such a category.

There are five sections. Section~\ref{sec:1} is aimed at someone who is completely new to the word ``motive''. In it, we discuss a proof of Manin showcasing how the category $\sM_{k, \QQ}^{\eff}$ of classical motives can be used to give a conceptual explanation of the numbers of rational points of certain varieties over finite fields, cf.Eq.\ref{equa:vn}, Eq.\ref{equa:manin2}. 

In Section~\ref{sec:2} we connect $\sM_{k, \QQ}^{\eff}$  to Voevodsky's motivic cohomology $H^i(X, \ZZ(j))$ for smooth schemes $X$ via Voevodsky's triangulated category of motives $\DM_k^{\eff}$. The most important parts of Section~\ref{sec:2} are the definition of the complex of presheaves $\ZZ(i)$ on $\Sm_k$, Eq.\eqref{equa:ZZiL}, and the Atiyah-Hirzebruch spectral sequence, Eq.\eqref{equa:HMKSS}, as these will be used later. 

In Section~\ref{sec:3} we discuss the cdh topology as one approach to working with motives of non-smooth varieties, and motives with compact support. Again, this is aimed at someone who has never seen the term cdh before. We present nilinvariance as a shortcoming of cdh motivic cohomology, Exam.\ref{exam:cdhBad}, and introduce the coarser \emph{procdh topology}, Def.\ref{defi:procdh}, as a solution.

In Section~\ref{sec:4} we describe Elmanto-Morrow's motivic cohomology, Def.\ref{defi:EM}, and present the equivalence $\ZZ(n)^{\procdh} \cong \ZZ(n)^{\EM}$, Thm.\ref{theo:procdhEM}. Both $\ZZ(n)^{\procdh}$ and $\ZZ(n)^{\EM}$ sit in Atiyah-Hirzebruch spectral sequences, Prop.\ref{prop:procdhAH}, Rem.\ref{rema:ZEMrem}.

In Section~\ref{sec:5} we mention recent work of Annala, Iwasa and Annala, Iwasa, Hoyois on a category of non-$\AA^1$-invariant motivic homotopy types where $\ZZ(n)^{\EM}$ is representable.

I thank the organisers of the workshop for the wonderful opportunity to visit RIMS in Kyoto and speak about this work. I also thank Federico Binda, Tess Bouis, Thomas Geisser, Annette Huber, Ryomei Iwasa, Matthew Morrow, Kay R{\"u}lling, Shuji Saito, and Vladimir Sosnilo for comments during the preparation of these expanded notes. I take full responsibility for any places where I have ignored or misunderstood their good advice and ended up introducing inaccuracies and mistakes.

\section{What is a motive, though?} \label{sec:1}

It is never clear how much background the audience has to the theory of motives. In an attempt to ground the talk in the ``real world'' we begin with a result from the origins of the theory. This doesn't have much to do with the main theorems of the talk, but hopefully will give the reader some kind of picture to associate with the words ``motive'' and ``motivic cohomology''. The reader not interested in classical motives can go directly to Section 2.

%
%

\begin{theo}[{Manin, 1968, \cite{Man68}}] \label{theo:manin1}
Suppose that $X$ is a smooth projective unirational threefold over a finite field $k = \FF_q$. Then the number $ν_n$ of rational points over $\FF_{q^n}$ is 
\begin{equation} \label{equa:manin}
ν_n = 1 + q^n - q^n\sum_{i \in I} ω_i^n  + q^{2n} + q^{3n}
\end{equation}
for some finite set $I$ where $ω_i$ are algebraic integers with $|ω_i| = q^{1/2}$.
\end{theo}

\begin{rema}
In the talk, Eq.\eqref{equa:manin} was incorrectly stated as $ν_n = 1 + q^n - q^n\sum_{i = 1}^{10} ω_i^n  + q^{2n} + q^{3n}$. The 10 in the sum appears in the special case where $X$ is a cubic hypersurface in $\PP^4$; it is twice the dimension of the Albanese variety $A(S)$ of the surface $S$ parametrising the lines of $X$. 
This special case was proved by Bombieri and Swinnerton-Dyer in \cite{BS67}, and is presented by Manin as an application of his more general result.
Indeed, Bombieri and Swinnerton-Dyer state their result as $ν_n = 1 + q^n - q^n Tr(π^n) + q^{2n} + q^{3n}$ where $Tr(π^n)$ is the trace of the $n$th power of the Frobenius endomorphism of $A(S) / \FF_q$ in the endomorphism algebra of $A(S)$ over $\QQ$.
\end{rema}



\begin{proof}[Proof. (Manin)]
Unirational, by definition, means there exists a generically finite morphism $U \to X$ from some dense open subscheme $U \subseteq \PP^d$. Of course $d = 3$ in our case. By results of Abhyankar on resolution of singularities of threefolds in positive characteristic, \cite[pg.1, ``Dominance'']{Abh66}, %
the morphism $U \to X$ can be promoted to a roof
\[ \PP^3 = Y_0 \stackrel{p_1}{\leftarrow} Y_1 \stackrel{p_2}{\leftarrow} \dots \stackrel{p_n}{\leftarrow} Y_n = Y \stackrel{f}{\to} X \]
where $f$ is generically finite and each $p_i$ is the blowup in a smooth connected centre.
Writing
$H^{even}(X) = \oplus_r H_{\et}^{2r}(X {\times_{\FF_q}} \overline{\FF}_q, \QQ_\ell)$%
, 
$H^{odd}(X) = \oplus_r H_{\et}^{2r+1}(X {\times_{\FF_q}} \overline{\FF}_q, \QQ_\ell)$%
, and $H^*(X) = H^{even}(X) \oplus H^{odd}(X)$,  
we recall that the Lefschetz trace formula, \cite[pg.86,Thm.3.2]{SGA45}, says 
\[ ν_n = 
\sum \left \{ \begin{array}{c} \textrm{eigenvalues of} \\ \operatorname{Frobenius}|{H^{even}(X)} \end{array} \right \}
 - 
\sum \left \{ \begin{array}{c} \textrm{eigenvalues of} \\ \operatorname{Frobenius}|{H^{odd}(X)} \end{array} \right \}
 . \]
Poincaré duality gives a Gysin morphism $f_*: H^*(Y) \to H^*(X)$ satisfying $f_*f^* = \deg f \cdot \id$, \cite[Exp.XVIII]{SGA43}. So 
\begin{align} \label{equa:summandFF}
H^*(Y) &= H^*(X) \oplus \{\textrm{ something }\}
\end{align}
in the category 
\[ 
\Gal_{\QQ_\ell} = 
\left \{ 
\begin{array}{c} 
\textrm{graded } \QQ_\ell \textrm{-vector spaces } \\ 
\textrm{ with a continuous } \hat{\ZZ}\textrm{-action} 
\end{array} 
\right \}
\]
Since $\dim Y_i = 3$, the centre $Z_i$ of the blowup $Y_{i+1} \to Y_i$ is a point or a curve. The blowup formula, \cite[0EW4]{stacks-project},
then says 
\begin{align} \label{equa:PBC}
H^*(Y_{i+1}) &= 
\left \{ 
\begin{array}{ll} 
H^*(Y_i) \oplus H^*(Z_i)(-1) & \dim Z_i = 1 \\
H^*(Y_i) \oplus H^*(Z_i)(-1) \oplus H^*(Z_i)(-2) & \dim Z_i = 0
\end{array} 
\right . 
\end{align}
where 
\[ -(-j) = - \otimes H^{2}(\PP^1)^{\otimes j} \]
Finally, we know that
\begin{align} \label{equa:Proj}
H^*(\PP^3) =& \oplus_{j = 0}^3 \QQ_\ell(-j).
\end{align}

By Eq.\eqref{equa:summandFF} the eigenvalues of $H^*(X)$ are among the eigenvalues of $H^*(Y)$. And using Eq.\eqref{equa:PBC},  Eq.\eqref{equa:Proj}, and induction on $i$ we reduce to calculating the eigenvalues of $H^*(W)$ for $\dim W = 0, 1$. Dimension zero is easy and dimension one was done by Weil in the 40's, 
\cite{We49}.
\end{proof}

\begin{obse}\ 

\begin{center}
\fbox{
\parbox{5cm}{
\begin{center}
The decompositions 
\eqref{equa:summandFF}, \eqref{equa:PBC}, \eqref{equa:Proj}
all come from algebraic cycles.
\end{center}
}
}
\end{center}
\end{obse}

Explicitly, consider the cycle class map $cl: \sZ^r(X) \to H^{2r}(X, \QQ_\ell(r))$ where $\sZ^r(X)$ is the free abelian group of closed irreducible subschemes $Z \subseteq X$ of codimension $r$.
Any cycle $γ \in \sZ^r(X {\times} Y)$ with $r = \dim Y$ induces a morphism of graded Galois modules 
\begin{align*}
H^*(γ): H^*(Y, \QQ_\ell) &\to H^{*}(X, \QQ_\ell) \\
	α &\mapsto pr_{1*}(cl(γ) \cap pr_2^*α)
\end{align*}
where $pr_1, pr_2: X \times Y \to X, Y$ are the two projection morphisms, and $\cap$ is the product on the bigraded abelian group $H^*_{\et}(X \times Y, \QQ_\ell(\star))$.
%
%

\begin{exam} \label{exam:cyclesOnEtale}\ 
\begin{enumerate}
 \item 
 In the case of Eq.\eqref{equa:summandFF}, the direct summand $H^*(X)$ of $H^*(Y)$ is the image of the morphism associated to $\tfrac{1}{\deg f} [Y \times_X Y] \in  \sZ^{\dim Y}(Y \times Y)$.


 \item In the case of Eq.\eqref{equa:Proj}, the direct summand $\QQ_\ell(j)$ is cut out by the cycle $[\PP^{3-j} {\times} \PP^{j}]$ where $\PP^k \subseteq \PP^n$ is the embedding $[(\ast, \dots, \ast, 0, \dots, 0)]$. 

\end{enumerate}
\end{exam}

\begin{defi}[Grothendieck, 60's]
Cf. \cite[\S 2, \S 3, \S 4]{Man68}. Define the $\QQ$-linear additive category $\mathcal{C}_{k, \QQ}$ as follows. Objects are smooth projective varieties. Given smooth projective varieties $Y, X$ we set $\hom_{\mathcal{C}_{k, \QQ}}(Y, X) = CH^{\dim Y}(X {\times} Y)_\QQ$, where $CH^d(-)_\QQ$ means the $\QQ$-vector space of cycles of codimension $d$ modulo rational equivalence, \cite[\S 1.3]{Ful84}.\footnote{More generally, if $Y$ is not equidimensional then we use $\oplus_{i} CH^{\dim Y_i}(X {\times} Y_i)_\QQ$ where the $Y_i$ are the connected components of $Y$.} Composition is, \cite[\S 16]{Ful84}, 
\begin{align*}
CH^{\dim X_2}(X_1 \times X_2)_\QQ \times CH^{\dim X_3}(X_2 \times X_3)_\QQ &\to CH^{\dim X_3}(X_1 \times X_3)_\QQ \\
(α, β) &\mapsto pr_{13*} (pr_{12}^*α \cdot pr_{23}^*β)
\end{align*} 
where $pr_{ij}: X_1 \times X_2 \times X_3 \to X_i \times X_j$ are the canonical projections, $(-)^*$ means flat pullback and $(-)_*$ is proper pushforward, and $\cdot$ is the intersection product, \cite[\S 1.4, \S 1.7, \S 6]{Ful84}.

Cf. \cite[\S 5]{Man68}. The category $\mathcal{C}_{k, \QQ}$ is an additive category where disjoint union of varieties corresponds to direct sum. However, idempotents in the category $\mathcal{C}_{k, \QQ}$ do not necessarily split in the sense that there are morphisms $ε: M \to M$ satisfying $ε^2 = ε$, but not admitting a decomposition $M {\substack{\danger \\ =}} M_1 \oplus M_2$ for which $ε$ is $\id_{M_1} \oplus\ 0_{M_2}$. For example, $\PP^3$ is connected, so considered as an object of $\mathcal{C}_{k, \QQ}$ it has no proper nontrivial direct summands. 

We formally impose such decompositions by passing to the idempotent completion $\mathcal{M}^{\eff}_{k, \QQ} = \mathcal{C}^\natural_{k, \QQ}$. Explicitly, an object of $\sM^{\eff}_{k, \QQ}$ is a pair 
\[ (X, ε) \]
where $X \in \SmProj_k$ and $ε \in CH^{\dim X}(X {\times} X)_\QQ$ satisfies $ε \circ ε = ε$. Heuristically, the pair $(X, ε)$ represents the direct summand of $X$ cut out by $ε$, so, unsurprisingly, the new hom sets are $\hom_{\mathcal{M}^{\eff}_{k, \QQ}}((X_1, ε_1), (X_2,ε_2))  = \{ f \in \hom_{\mathcal{C}_{k, \QQ}}(X_1, X_2)\ |\ ε_2 \circ f \circ ε_1 = f\}$. Then for any such $(X, ε)$ we have the canonical decomposition $(X, \id) \cong (X, ε) \oplus (X, 1{-}ε)$.%
\footnote{More generally, given an idempotent $ε = η \circ ε\circ η: (X, η) \to (X, η)$ we have $(X, η) \cong (X, η) \oplus (X, η{-}ε)$. }%

Sending a morphism $f: X \to Y$ to the cycle associated to its graph\footnote{Explicitly, $\Gamma_f = im(\id \times f: X \to X \times Y)$.} $[\Gamma_f] \in CH^{\dim Y}(X {\times} Y)_\QQ$ produces a \emph{contra}variant functor
\[ \SmProj_k^{\op} \to \sM_{k, \QQ}^{\eff} \]
from the category of smooth projective $k$-varieties. For this section we will use the same notation for a variety and its image in $\sM_{k, \QQ}^{\eff}$ since it will be clear from the context which is intended.

\end{defi}

\begin{exam}
Since the diagonal $[diag.] \in CH^1(\PP^1 {\times} \PP^1)$ is rationally equivalent to $p_0 + p_1$ where $p_0 = [\PP^1{\times}\{0\}]$ and $p_1 = [\{0\}{\times}\PP^1]$, and $p_0$, $p_1$ are idempotents, we obtain a decomposition $\PP^1 \cong (\PP^1, p_0) \oplus (\PP^1, p_1)$ in $\sM^{\eff}_{k, \QQ}$. The former summand is isomorphic to $\QQ := \Spec(k)$. The latter summand is denoted $\LL$ and called the \emph{Lefschetz motive}. That is,
\[ \PP^1 \cong \QQ \oplus \LL. \]
\end{exam}

\begin{exam}
More generally, for any irreducible smooth projective variety $X$ of dimension $d$ and any closed point $x \in X$ we have two canonical idempotents, $p_0 = \tfrac{1}{[k(x):k]} [X {\times} x]$ and $p_d = \tfrac{1}{[k(x):k]} [x {\times} X]$ in $CH^d(X {\times} X)_\QQ$. If $d > 0$ these are orthogonal in the sense that $p_0 \circ p_d = p_d \circ p_0 = 0$. Correspondingly, we obtain a decomposition $X \cong (X, p_0) \oplus (X, \id - p_0 - p_d) \oplus (X, p_d)$ in $\sM_{k, \QQ}^{\eff}$. The left and right summands are isomorphic to $\QQ$ and $\LL^{\otimes d}$ respectively, the tensor product on $\sM^{\eff}_{k, \QQ}$ being defined by $(X, ε) \otimes (Y, η) = (X \times Y, ε \boxtimes η)$ where $\boxtimes: CH^d(X {\times} X)_\QQ {\times} CH^e(Y {\times} Y)_\QQ \to CH^{d+e}(X {\times} Y {\times} X {\times} Y)_\QQ$ has the obvious meaning. The middle summand, $(X, \id - p_0 - p_d) $, is denote $X^+$ in \cite[\S 10]{Man68}. So we have
\[ X \cong \QQ \oplus X^+ \oplus \LL^{\otimes d}. \]
\end{exam}

This setup complete, we can now lift the decompositions 
\begin{align*}
H^*(Y) &= H^*(X) \oplus \{\textrm{ something }\} \\
H^*(Y_i) &= 
\left \{ 
\begin{array}{l} 
H^*(Y_i) \oplus H^*(Z_i)(-1)  \qquad \textrm{ or } \\
H^*(Y_i) \oplus H^*(Z_i)(-1) \oplus H^*(Z_i)(-2)  
\end{array}
\right . \\
H^*(\PP^3) &= \QQ_\ell \oplus \QQ_\ell(-1) \oplus \QQ_\ell(-2) \oplus \QQ_\ell(-3).
\end{align*}
in $\Gal_{\QQ_\ell}$ to decompositions in $\sM_{k, \QQ}^{\eff}$
\begin{align}
Y &= X \oplus 	\{\textrm{ something } \} \label{MYX} \\
Y_{i+1} &= 
\left \{ 
\begin{array}{l} 
Y_i\  \oplus \  Z_i {\otimes} \LL \quad \textrm{ or } \\
Y_i\  \oplus\  Z_i{\otimes} \LL\  \oplus\  Z_i{\otimes} \LL^{\otimes 2} 
\end{array} 
\right . \label{MYZ} \\
\PP^3 &= \QQ\ \oplus\ \LL\ \oplus\ \LL^{\otimes 2}\ \oplus\ \LL^{\otimes 3}. \label{MPL}
\end{align}
Note the above is now independent of $\ell$! Applying the additive monoïdal functor 
\begin{align*}
H^*: \mathcal{M}^{\eff}_{k, \QQ} &\to Gal_{\QQ_\ell} \\
(X, ε) &\mapsto image \biggl ( H^*(ε): H^*_{et}(X, \QQ_\ell) \stackrel{}{\to} H^*_{\et}(X, \QQ_\ell) \biggr ) 
\end{align*}
recovers \eqref{equa:summandFF}, \eqref{equa:PBC}, \eqref{equa:Proj}. The expression 
\begin{equation} \label{equa:vn}
ν_n = 1 + q^n - q^n\sum_{i \in I} ω_i^n  + q^{2n} + q^{3n}
\end{equation}
in Theorem~\ref{theo:manin1} now becomes a consequence of the following theorem, combined with Weil's proof of the Riemann Hypothesis for curves, \cite{We49}.

\begin{theo}[{\cite[\S 11]{Man68}}] \label{theo:ManMot}
Suppose that $X$ is a smooth projective unirational three-fold over a finite field $k = \FF_q$. Then there is a decomposition
\begin{equation} \label{equa:manin2}
X\ \cong\ 
\QQ \ \oplus\ \LL\ \oplus\ \LL{\otimes}U\ \oplus\ \LL^{\otimes 2}\ \oplus\ \LL^{\otimes 3} 
\end{equation}
in $\mathcal{M}^{\eff}_{k, \QQ}$ where $U$ is a direct summand of a motive of the form $\oplus_{j \in J} C_j^+$, the $C_j$ being curves.
\end{theo}

\begin{rema}
It turns out that the decompositions Eq.\eqref{MYX}, Eq.\eqref{MYZ}, Eq.\eqref{MPL} work over any base field. If the base is some subfield $K \subseteq \CC$, then sending $X$ to the sum of its de Rham cohomology defines a functor 
\begin{align*}
\sM_{k, \QQ}^{\eff} &\to \Hdg := \{ \textrm{ graded pure Hodge structures } \} \\
X &\mapsto H^*_{\dR}(X)
\end{align*}
towards the category of graded pure Hodge structures. Then the decomposition of Theorem~\ref{theo:ManMot} has consequences, for example, for the periods of $X$.

If the base is a field of positive characteristic $p$, then crystalline cohomology gives a functor towards the category of graded $W(k)[\frac{1}{p}]$-modules equipped with a semi-linear endomorphism.
\begin{align*}
\cM^{\eff}_{k, \QQ} &\to \Crys := \left \{ \begin{array}{c} W(k)[\tfrac{1}{p}]\textrm{-modules equipped with } \\ \textrm{ a semi-linear endomorphism } \end{array} \right \}	 \\
X &\mapsto H_{\crys}^*(X)
\end{align*}
and similarly, the decompositions Eq.\eqref{MYX}, Eq.\eqref{MYZ}, Eq.\eqref{MPL} have implications for $H_{\crys}^*(X)$.
\end{rema}

Here is a picture:
\begin{equation} \label{equa:reas}
\xymatrix@R=2ex{
&& \Gal_{\QQ_\ell} & \FF_\ell \not\subseteq k \\
\SmProj_k^{\op} \ar[r] & \cM^{\eff}_{k, \QQ} \ar[ur]^{H^*_{\ell{\operatorname{-adic}}}} \ar[r]^(0.6){H^*_{\dR}} \ar[dr]_{H^*_{\crys}} & \Hdg & \QQ \subseteq k \subseteq \CC \\
&& \Crys & \FF_p \subseteq k
}
\end{equation}

\section{Smooth motives and $K$-theory} \label{sec:2}

In this section we discuss Voevodsky motives. The reader familiar with this theory can jump to Section~\ref{sec:3}. The most important part of this section is the definition of the complex of presheaves $\ZZ(i)$ on $\Sm_k$, Eq.\eqref{equa:ZZiL}, and the Atiyah-Hirzebruch spectral sequence, Eq.\eqref{equa:HMKSS}, as these will be used later.

In the 90's, Voevodsky defined a \emph{co}variant functor $M: \Sm_k \to \DM^{\eff}_k$ towards a triangulated category fitting into a commutative square, \cite[Prop.2.1.4]{Voe00},
\[ \xymatrix{
\SmProj_k \ar[d] \ar[r] & (\cM^{\eff}_k)^{\op} \ar[d]^{\textrm{fully faithful}} & \\
\Sm_k \ar[r]_M & \DM^{\eff}_k & 
} \]
The functor $M$ is universal (in a sense that can be made mathematically precise) with respect to the properties:
\begin{enumerate}
 \item (Colimits) $\DM_k^{\eff}$ admits all small sums. In particular, idempotents split.%
\footnote{Explicitly, if $ε: X \to X$ satisfies $ε^2 = ε$, then $X \cong X_ε \oplus X_{1-ε}$ where $X_η = Cone(\oplus_\NN X \stackrel{η - σ}{\to} \oplus_\NN X)$ and $η - σ$ is the obvious interpretation of $(x_0, x_1, x_2, \dots) \mapsto (ηx_0, ηx_1 - x_0, ηx_2 - x_1, \dots)$.} %
 \item (Transfers) For certain finite morphisms $f: Y \to X$ there are ``backwards'' morphisms 
\[ {}^tf: M(X) \to M(Y). \]

 \item (Homotopy invariance) The morphisms 
 \[ M(\AA^1_X) \to M(X) \]
 are isomorphisms.
 
 \item (Zariski descent) For open immersions $U, V \subseteq X$ there is a distinguished triangle 
 \[ M(U \cap V) \to M(U) \oplus M(V) \to M(U \cup V) \to M(U \cap V)[1]. \]
\end{enumerate}

\begin{defi}[{\cite[pg.20]{Voe00}}] \label{defi:DM}
Concretely, to build $\DM_k^{\eff}$ Voevodsky first defines a mild enlargement $\Sm_k \subseteq \SmCor_k$ %
of the category of smooth $k$-varieties. 
Then 
\[ \DM_k^{\eff} \subseteq D(\PSh(\SmCor_k, \ZZ)) \]
is the smallest full sub-triangulated category admitting all small sums and containing the complexes of presheaves on $\SmCor_k$
\begin{equation} \label{eq:MXsing}
M(X) := L^{\Zar}L^{\AA^1}\hom_{\SmCor_k}(-, X)
\end{equation}
for $X \in \Sm_k$, where $L^{\AA^1}$ means localisation with respect to the morphisms $\AA^1_T {\to} T$ for $T \in \Sm_k$ and $L^{\Zar}$ means localisation with respect to Zariski hypercoverings.\footnote{Or equivalently, since the Zariski topos is hypercomplete, localisation with respect to Zariski \v{C}ech hypercoverings.}
\end{defi}

\begin{rema}\ 
\begin{enumerate}
 \item In the top right corner we have used $\sM_k^{\eff}$ in place of the $\sM_{k, \QQ}^{\eff}$ that appeared in Section~\ref{sec:1}. It is defined in the same way, except one uses Chow groups with integral coefficients $CH^*(-)$ instead of rational coefficients $CH^*(-)_\QQ$. We used rational coefficients in Section~\ref{sec:1} to get the splitting Eq.\eqref{MYX}. The splittings Eq.\eqref{MYZ} and Eq.\eqref{MPL} hold integrally.


 \item Traditionally, $L^{\AA^1}K(-) = K(\Delta^\bullet \times -)$ is taken for $L^{\AA^1}$ where $\Delta^n = \Spec(\ZZ[t_0, \dots, t_n] / 1 = \sum t_i)$ and the differentials come from the alternating sum of the face morphisms $\Delta^n \to \Delta^{n+1}$, as in singular cohomology.

 \item The functor $L^{\Zar}$ is Zariski hypercohomology in the sense that for any $K \in D(\PSh(\SmCor_k, \ZZ))$ in the image of $L^{\AA^1}$, if one restricts $L^{\Zar}K$ to the small Zariski site $T_{\Zar}$ of  any $T \in \Sm_k$, one gets the Zariski hypercohomology $L^{\Zar}(K)|_{T_{\Zar}} = \HH_{\Zar}^*(K|_{T_{\Zar}})$, \cite[Prop.3.1.8, Thm.3.1.11]{Voe00}.
\end{enumerate}
\end{rema}

As in the case of pure motives, when $\FF_\ell \not\subseteq k$, for each $n \in \NN$ we have an étale realisation functor
\begin{align*}
(\DM_k^{\eff})^{\op} &\to D(k_{\et}, \ZZ / \ell^n) \\
M(X) & \mapsto Rf_*(\ZZ / \ell^n)
\end{align*}
where $f: X \to \Spec(k)$ is the structural morphism of $X$ and $\ZZ / \ell^n \in D(X_{\et}, \ZZ / \ell^n)$ is the constant sheaf, and when $\QQ \subseteq k \subseteq \CC$ we have a de Rham realisation functor
\begin{align*}
(\DM_k^{\eff})^{\op} &\to D(k) \\
M(X) &\mapsto R\Gamma_{\Zar}(X, \Omega^*_X)
\end{align*}
which can be embellished to capture mixed Hodge structures in the sense of \cite{Del74}.
References for these two realisations include \cite{Hub00}, \cite{Hub04}, \cite{Ivo07}, \cite{Lec08}, \cite{LW09}, \cite{Sch12}, \cite{LW13},\cite{Har16}, \cite[App.A]{KS17}.

The $p$-adic picture is more subtle. Crystalline cohomology is not $\AA^1$-invariant so by definition cannot give a functor out of $\DM_k^{\eff}$. Berthelot's rigid cohomology does give a functor by passing through étale motives without transfers, \cite[App.B]{Ayo14},  %
%
%
%
%
but this has $K$-coefficients where $K = W(k)[\tfrac{1}{p}]$, so the all-important $p$-torsion is missing. %
However, recently, integral $p$-adic realisation functors from $\DM_k^{\eff}$ have been constructed by Merici, \cite{Mer22}, and Annala, Hoyois, Iwasa, \cite{AHI24}. Under resolution of singularities, Merici's sends $M(X)$ to the log crystalline cohomology of a smooth compactification $X \subseteq \ol{X}$, and it is expected that the two agree.
\begin{align*}
(\DM_k^{\eff})^{\op} &\to D(W(k)) \\
M(X) &\mapsto R\Gamma_{\crys}((\overline{X}, \partial X)/W(k))
\end{align*}
Here $\partial X = \overline{X} \setminus X$. So now we (almost) have the picture:
\begin{equation} \label{equa:DMreal}
\xymatrix@R=2ex{
&& D(k_{\et}, \ZZ / \ell^n) & \FF_\ell \not\subseteq k \\
\Sm_k^{\op} \ar[r]^-M & (\DM_k^{\eff})^{\op} \ar[ur]^{R\Gamma_{\et}} \ar[r]^(0.6){R\Gamma_{\dR}} \ar@<0ex>[dr]^{R\Gamma_{\textrm{AHI}}} \ar@<-0.5ex>[dr]_{R\Gamma_{\textrm{Mer,}}} 
& D(k) & \QQ \subseteq k \subseteq \CC \\
&& D(W(k)) & \FF_p \subseteq k
}
\end{equation}
As in $\sM_{k,\QQ}^{\op}$ we have liftings to $\DM_k^{\eff}$ of various properties of these cohomology theories. One which will be important in this talk is a distinguished triangle associated to a blowup\footnote{Since we are assuming $Z$ and $X$ are smooth, this triangle is split, but we don't use this here.}
\begin{equation} \label{equa:DMbu}
M(E) \to M(Z) \oplus M(Bl_X Z) \to M(X) \to M(E)[1]
\end{equation}
where $Z {\to} X$ is a closed immersion in $\Sm_k$ and $E = Z {\times_X} Bl_X Z$, \cite[Prop.3.5.3]{Voe00}.

\begin{rema} \label{rema:missing}
Since $M(X) \cong M(X {\times} \AA^1)$ in $\DM^{\eff}_k$, we are missing other functors that we might like such as
\[ H_{\et}^*(-, \ZZ/p): (\DM^{\eff}_{\FF_p})^{\op} \stackrel{\danger}{\not\to} D(\FF_p), 
\qquad \qquad 
H_{\Zar}^*(-, \OO): (\DM^{\eff}_{\QQ})^{\op} \stackrel{\danger}{\not\to} D(\QQ) \]
We will come back to this in Section~\ref{sec:4}.
\end{rema}

We finish this section with a short discussion about $K$-theory. In the category $\sM_{k, \QQ}^{\op}$ for an irreducible smooth projective variety $X$ we have, almost by definition,
$ \hom_{\sM_{k}^{\op}}(X, \LL^{\otimes i}) \cong CH^i(X)$.
In $\DM_k^{\eff}$ this generalises to
\[ \HH_{\Zar}^i(X, \ZZ(j)) \cong \hom_{\DM_k^{\eff}}(M(X), \ZZ(i)[j]) \cong CH^i(X, 2i{-}j), \]
\cite{Voe02}, where $CH^*(X, \star)$ are Bloch's higher Chow groups, \cite{Blo86}, \cite{Lev94}. Here, one defines%
\begin{equation} \label{equa:ZZiL}
\ZZ(i) := \LL^{\otimes i}[-2i]
\end{equation}
where $\LL^{\otimes i}$ is considered as an object of $\DM_k^{\eff}$ under the embedding $(\sM_k^{\eff})^{\op} \subseteq \DM_k^{\eff}$. Explicitly, $\ZZ(i)$ is the direct summand of $L^{\Zar}L^{\AA^1}\hom_{\SmCor_k}(-, (\PP^1)^{\times i})$ complementary to all the ``axes'' $\PP^1 \times \dots \times \PP^1 \times \{\infty\} \times \PP \times \dots \times \PP^1$.%

Even better, for $X$ smooth, the Grothendieck-Riemann-Roch isomorphism $K_0(X) \otimes \QQ \cong CH^*(X) \otimes \QQ$, %
\cite{SGA6}, gets promoted to an Atiyah-Hirzebruch style spectral sequence, \cite{FS02}, \cite{Lev08}, originally conjectured by Beilinson, \cite{Bei82}, \cite{Bei87},
\begin{equation} \label{equa:HMKSS}
E_2^{i,j} := \HH_{\Zar}^{i-j}(X, \ZZ(-j)) \implies K_{-i-j}(X).
\end{equation}
Here $\ZZ(-j) = 0$ for $j > 0$. In particular, $E_{2}^{i,j} = 0$ for $j > 0$. It is also known that $E_2^{i,j} = 0$ for $i > \dim X$, and $i > -j$.

 It is this spectral sequence which motivates our definition of motivic cohomology for non-smooth schemes below.

\section{Singular motives and the cdh-topology} \label{sec:3}

In this section we discuss the cdh topology, followed by the procdh topology. The reader familiar with the cdh topology can skip to Definition~\ref{defi:procdh}.

So much for smooth varieties. What about singular varieties? The presheaves $\hom_{\Cor_k}(-, X)$ defined as 
\[ T \mapsto \ZZ \left \{Z \subseteq T \times X\ \middle |\ \begin{array}{c} Z \textrm{ is closed, irreducible, } \\ \textrm{ and finite surjective over an } \\ \textrm{  irreducible component of } X \end{array} \right \} \]
make sense more generally for any $X$ of finite type over $k$, however not much can be said of the corresponding objects $L^{\Zar}L^{\AA^1}\hom_{\Cor_k}(-, X)$ in $\DM_k^{\eff}$ without assuming that the base field satisfies some kind of resolution of singularities.\footnote{For resolution of singularities in the classical sense, this is done in \cite[\S 4]{Voe00}. In the author's PhD thesis, \cite{Kel17}, using Gabber's work on alterations, all the results from \cite[\S 4]{Voe00} are obtained in positive characteristic at the expense of inverting $p$ in the coefficients.} If the base field \emph{does} satisfy resolution of singularities, then we obtain an equivalence, \cite[Thm.4.1.2]{Voe00},
\[ \DM_k^{\eff} \cong \DM_{k, \cdh}^{\eff} \]
with the category of \emph{cdh-motives}, implicit in Voevodsky's work, cf.\cite{Voe10}, and explicitly studied by Cisinski and Déglise, \cite{CD15}. To motivate this topology let's talk about cohomology with compact support following the author's note \cite{Kel17}. 

In general, for any closed immersion $Z \to X$ with open complement $U \subseteq X$, a reasonable cohomology theory with compact support can be expected to have a localisation distinguished triangle of the form
\begin{equation} \label{equa:RURX}
 R\Gamma_c(U) \to R\Gamma_c(X) \to R\Gamma_c(Z) \to R\Gamma_c(U)[1].
 \end{equation}
Another property that is reasonable to expect is that for proper schemes $X$ (or compact topological spaces) the compact support version agrees with usual cohomology
\begin{equation} \label{equa:RXRX}
 R\Gamma_c(X) \cong R\Gamma(X), \qquad X \textrm{ proper.}	
\end{equation}
Given that we would like Eq.\eqref{equa:RURX} and Eq.\eqref{equa:RXRX}, a reasonable way to define $R\Gamma_c(X)$ starting from a functor $R\Gamma \in D(\PSh(\Sch_k, \ZZ))$ on the category $\Sch_k$ of separated schemes of finite type over $k$, is to choose a compactification $X \subseteq \ol{X}$ and then define
\[ R\Gamma_c(X) := \Cone \biggl ( R\Gamma(\ol{X}) \to R\Gamma(\partial X) \biggr )[-1]. \]
An obvious problem is that this is only well-defined if it is independent of the choice of compactification.
\begin{equation} \label{equa:coneIso}
\Cone \biggl (R\Gamma(\ol{X}) \to R\Gamma(\partial X) \biggr ) \stackrel{?}{\cong} \Cone \biggl ( R\Gamma(\ol{X}') \to R\Gamma(\partial X') \biggr ).
\end{equation}
For this question of independence of the compactification, it suffices to consider the case that there is a morphism $f: \ol{X}' \to \ol{X}$ compatible with the embeddings $X \subseteq \ol{X}', X \subseteq \ol{X}$. 
In this case, the isomorphism Eq.\eqref{equa:coneIso} is equivalent to asking that the sequence
\[ R\Gamma(\ol{X}) \to R\Gamma(\ol{X}') \oplus R\Gamma(\partial X) \to R\Gamma(\partial X') \]
fits into a distinguished triangle.%
\footnote{Of course, at the level of triangulated categories none of this makes sense unless everything is rigidified in some way. We should instead be talking about fibre and cofibre sequences in some kind of stable $\infty$-categories, for example pretriangulated dg-categories.} %
Noting that this resembles a Mayer-Vietoris condition, one could try and turn this into a sheaf condition, and this is precisely what Voevodsky does.

\begin{defi}[Voevodsky] \label{defi:cdhDef}
Write $\Sch_k$ for the category of separated schemes of finite type over $k$. The \emph{cdh topology} on $\Sch_k$ is generated by families of the form
\begin{enumerate}
 \item Nisnevich coverings.
 \item \label{abu} Families of the form 
 \[ \{Z \to X, Y \to X\} \]
 where $Z \to X$ is a closed immersion, $p: Y \to X$ is proper, and $p$ is an isomorphism outside of $Z$, i.e., $p^{-1}(X \setminus Z) \stackrel{\sim}{\to} X \setminus Z$. 
\end{enumerate}
\end{defi}

\begin{prop}[{\cite{Voe10}}] \label{prop:MVcdh}
A complex of presheaves $C \in D(\Shv_{\Nis}(\Sch_k, \ZZ))$ lies in $D(\Shv_{\cdh}(\Sch_k, \ZZ))$ if and only if for every $Z, Y, X$ as in Def.\ref{defi:cdhDef}\eqref{abu}, the sequence
\[ C(X) \to C(Y) \oplus C(Z) \to C(Z \times_X Y) \]
fits into a distinguished triangle.
\end{prop}

\begin{rema}
Following the discussion above, another way to phrase Proposition~\ref{prop:MVcdh} would be: $C$ is in $D(\Shv_{\cdh}(\Sch_k, \ZZ))$ if and only if it is in $D(\Shv_{\Nis}(\Sch_k, \ZZ))$ and has a well-defined theory of compact support.
\end{rema}

It is a little more subtle, but using the material in \cite{SVRelCyc} one can make a version $\Cor_k$ of $\SmCor_k$ whose objects are the objects of $\Sch_k$ which contains $\SmCor_k$ as a full subcategory in the obvious way.\footnote{In the notation of \cite{SVRelCyc} one takes $\hom_{\Cor_k}(X, Y) = c_{equi}(X \times Y / X, 0)$.} One can find expositions of this in \cite{Ivo07}, \cite{CD19}, \cite{Kel17}, and now also the Stacks Project \cite{stacks-project}. Then one can build a category of cdh motives as follows.

\begin{defi} \label{defi:cdhMotives}
Take the smallest full sub-triangulated category
\[ \DM_{k, \cdh}^{\eff} \subseteq D(\PSh(\Cor_k, \ZZ)) \]
admitting all small sums and containing the complexes of presheaves%
\footnote{Strictly, speaking, $L^{\cdh}L^{\AA^1}$ is only known to be $\AA^1$-invariant under some kind of resolution of singularities, \cite[Thm.5.5(3)]{FV00}, \cite[Cor.5.2.3]{Kel17}. So in general one should take the localisation $L^{\cdh, \AA^1}$ with respect to both classes at once, rather than $L^{\cdh}L^{\AA^1}$.} %
 on $\Cor_k$
\begin{equation} \label{eq:MXsingSing}
M(X) = L^{\cdh}L^{\AA^1}\hom_{\Cor_k}(-, X).
\end{equation}
\end{defi}

As discussed above, this category of cdh-motives will have a good theory of motives with compact support, almost by definition. Some cohomology theories such as $\ell$-adic cohomology have cdh-descent and as such still admit realisation functors, but some do not.

\begin{exam} \label{exam:cdhBad}
In Definition~\ref{defi:cdhDef}\eqref{abu} one can take $Z = X_{\red}$. Then $Y = \varnothing$ and the fibre sequence condition of Proposition~\ref{prop:MVcdh} becomes $C(X) \cong C(X_{\red})$. However, algebraic $K$-theory is not nilinvariant in general. For example \[ K_1(k) = k^* \neq (k[ε]/ε^2)^* = K_1(k[ε] / ε^2). \]
 Hence, we cannot hope to have a spectral sequence such as Eq.\ref{equa:HMKSS} for general schemes in $\Sch_k$ using the motivic cohomology represented in $\DM_{k,\cdh}^{\eff}$.
\end{exam}

\begin{exam}
By work of Morrow, \cite{Mor16}, and Kerz-Strunk-Tamme, \cite{KST18}, algebraic $K$-theory \emph{does} have fibre sequences for blowup squares if we remember the formal neighbourhoods of the closed immersions. That is, suppose $Z, X, Y$ are as in Def.\ref{defi:cdhDef}\eqref{abu}. Define $Z_n = \underline{\Spec} (\OO_{X} / \sI_Z^n)$ as the $n$th infinitesimal neighbourhood of $Z$, and similarly, $E_n := Z_n \times_X Y$. Then there is a distinguished triangle of spectra
\[ 
K(X)
\to 
R\lim K(Z_n)
\oplus K(Y)
\to 
 R\lim K(E_n)
\to K(X)[1].
\]
\end{exam}

We make the following definition.

\begin{defi}[K., Saito, {\cite[Def.1.1]{KS24}}] \label{defi:procdh}
The \emph{procdh topology} on $\Sch_k$ is generated by families of the form:
\begin{enumerate}
 \item Nisnevich coverings.
 \item \label{proabu} Families of the form 
 \[ \{Z_n \to X\}_{n \in \NN} \sqcup \{Y \to X\} \]
 where $Z \to X$ is a closed immersion, $p: Y \to X$ is proper, $p$ is an isomorphism outside of $Z$, i.e., $p^{-1}(X \setminus Z) \stackrel{\sim}{\to} X \setminus Z$, and $Z_n = \underline{\Spec} (\OO_{X} / \sI_Z^n)$.
\end{enumerate}
\end{defi}

\begin{prop}[K., Saito, {\cite[Thm.1.2]{KS24}}] \label{prop:MVprocdh}
A complex of presheaves $C \in D(\PSh(\Sch_k, \ZZ))$ lies in $D(\Shv_{\procdh}(\Sch_k, \ZZ))$ if and only if it has Nisnevich descent, and for every $Y, X, Z_0, Z_1, \dots$ as in Def.\ref{defi:procdh}\eqref{proabu}, the sequence
\[ C(X) \to C(Y) \oplus R\lim C(Z_n) \to R\lim C(Z_n \times_X Y) \]
is a fibre sequence.
\end{prop}

\begin{rema} \label{rema:cdhA1}
Looking at the definitions, one sees that the cdh topology is generated by the procdh topology and families of the form $\{Y_{\red} \to Y\}$. On the other hand, by the definition of relative cycles, we have isomorphisms $Y_{\red} \cong Y$ in $\Cor_k$ for every $Y \in \Sch_k$. Consequently, for objects coming from $\PSh(\Cor_k)$, the cdh and procdh sheafifications will agree. This leads to an equivalence
\[ \DM_{k, \procdh}^{\eff} \cong \DM_{k,\cdh}^{\eff} \]
if the left hand side is defined as in Definition~\ref{defi:cdhMotives} with cdh replaced with procdh.
\end{rema}

Instead of using the complexes coming from motives of the form \eqref{eq:MXsingSing}, we start with the left Kan extension of the complexes $\ZZ(n)$ on $\Sm_k$.

\begin{defi} \label{defi:LKan}
Write $L^{\Kan}$ for the left adjoint of the adjunction
\[ L^{\Kan}: D(\PSh(\Sm_k)) \rightleftarrows D(\PSh(\Sch_k)): \textrm{restriction} \]
where the right adjoint is restriction.
\end{defi}

\begin{rema}
The functor $L^{\Kan}$ above is a more general version of a functor that the reader is perhaps more familiar with. Considers the analogous adjunction 
\[ L: \PSh(\Et_X) \rightleftarrows \PSh(\Sch_X):R \]
where $\Et_X$ is the category of étale schemes over an affine scheme $X$ and $\Sch_X$ the category of finite type $X$-schemes, the henselisation of $X$ along a closed subscheme $Z \subseteq X$ is 
\[ X^h_Z = \underline{\Spec} (L\OO(Z)). \]
That is, global sections of $X^h_Z$ are calculated as the left Kan extension of the structure sheaf $\OO: Y \mapsto \Gamma(Y, \OO_Y)$ evaluated at the closed immersion $Z \to X$. So in this way, the functor $L^{\Kan}$ from Definition~\ref{defi:LKan} can be considered as a kind of very general henselisation.

There is at least one major difference though, since $\Sm_k$ doesn't admit fibre products, the \mbox{colimits} calculating $L^{\Kan}$ are not filtered, so they need to be derived, cf.\cite{BK72}. To concretely describe $L^{\Kan}C(Y) = L\varinjlim_{\substack{Y \to X \\ X \in \Sm_k}}C(X)$ for some $C \in D(\PSh(\Sm_k))$ and $Y \in \Sch_k$ one can use an appropriate simplicial object as people do when working with the cotangent complex, cf.\cite[Chap.XII,\S5]{BK72}.
\end{rema}

%
%
%
%
%

\newcommand{\colim}{\operatorname{colim}}

\begin{defi}[K., Saito, {\cite[Def.1.9]{KS24}}] \label{defi:Znprocdh}
Define $\ZZ(n)^{\procdh}$ as 
\[ \ZZ(n)^{\procdh} := L^{\procdh}L^{\Kan}\ZZ(n) \]
where $L^{\procdh}$ is the left adjoint to the inclusion $D(\Shv_{\procdh}(\Sch_k)) \subseteq D(\PSh(\Sch_k))$. 
That is, the procdh sheafifiaction of the left Kan extension of $\ZZ(n) \in D(\PSh(\Sm_k))$ from smooth schemes to schemes of finite type. 
\end{defi}

By Bhatt-Lurie, connective $K$-theory for $k$-algebras is left Kan extended from smooth algebras, see \cite[Ex.A.0.6]{EHKSY}. Moreover, non-connective $K$-theory is the procdh sheafification of connective $K$-theory, \cite[Thm.1.8]{KS24}. It follows that if we apply $L^{\procdh}L^{\Kan}$ to any of the filtered presheaves that give rise to the Atiyah-Hirzebruch style spectral sequence Eq.\eqref{equa:HMKSS} we obtain a spectral sequence from procdh motivic cohomology to algebraic $K$-theory.
 
\begin{prop}[K., Saito, {\cite[Thm.1.10]{KS24}}] \label{prop:procdhAH}
For any $Y \in \Sch_k$ there exists a bounded spectral sequence
\begin{equation} \label{equa:procdhAHSS}
\HH_{\procdh}^{i-j}(Y, \ZZ(-j)^{\procdh}) \implies K_{-i-j}(Y)	
\end{equation}
converging to algebraic $K$-theory.
\end{prop}

\section{Trace methods} \label{sec:4}

So now we have one new motivic complex $\ZZ(n)^\procdh$ what about the other one indicated in the title. Recall from Remark~\ref{rema:missing} that we are missing realisation functors $H_{\et}^*(-, \ZZ/p)$ for $k = \FF_p$ and $H_{\Zar}^*(-, \OO)$ for $k = \QQ$ because these are not $\AA^1$-invariant. More importantly for our present purposes 
there are no functors out of $\DM_k^{\eff}$ whose composition with $M: \Sm_k^{\op} \to \DM_k^{\eff}$ gives 
\begin{align*}
R\Gamma_{\Zar}(-, \Omega^{\geq n}) & \textrm{ for } k = \QQ, \textrm{ or }\\
R\Gamma_{\et}(-,ν_r^n) & \textrm{ for } k = \FF_p.
\end{align*}
Here,
\[ ν_r^n = \Gamma_{\et}(-, im(\dlog: (\OO^*)^{\times n} \to W_r\Omega^n)) \]
is the étale (or Zariski) sheafification of the image of the map presheaves of abelian groups $(x_1, \dots, x_n) \mapsto \tfrac{d[x_1]}{[x_1]}  \wedge \dots \wedge \tfrac{d[x_n]}{[x_n]}$. Here $[-]: \OO \to W_r\OO$ denotes the Teichm{\"u}ller map, cf.\cite[\S1.2, \S 2.2, Cor.4.2(i)]{Mor19}.

We do however have canonical comparison maps 
\begin{align}
\ZZ(n) &\to \Omega^{\geq n} \label{equa:ZZnOmega} \\	
\ZZ(n) &\to \lim_r ν_r^n[-n] \label{equa:ZZnnu}
\end{align}
on $\Sm_k$ coming from the canonical element $\tfrac{dt_1}{t_1} \wedge \dots \wedge \tfrac{dt_n}{t_n}$ in $H_{\dR}^n(\GG_m^{\times n})$ which we describe now.

\begin{cons} \label{cons:ZZnOmega}
Given a ring%
\footnote{Usually Milnor $K$-theory is only considered for fields, or local rings with infinite residue field, since it doesn't have the desired properties for a general ring. We are working over $\QQ$ and about to take the Zariski sheafification though.} %
 $A$, write $K_*^M(A)$ for the quotient of the tensor $\ZZ$-algebra $\oplus_{n \in \NN} (A^*)^{\otimes n}$ by the ideal generated by elements of the form $a \otimes (1-a)$ for $a \in A \setminus \{0, 1\}$, \cite{Mil70}. This is functorial in $A$, and we write $\sK_n^M$ for the Zariski sheaf associated to the presheaf $X \mapsto K_n^M(\Gamma(X, \OO_X))$ on $\Sm_\QQ$. Sending $a \in A^*$ to $\tfrac{da}{a} \in \Omega^1_A$ defines a morphism of graded rings $K_*^M(A) \to \Omega^*_A$ since $\tfrac{da}{a} \wedge \frac{d(1-a)}{1-a} = 0$. This induces a morphism of Zariski sheaves $\sK_n^M \to \Omega^n$ on $\Sm_\QQ$. Since the cohomology Zariski sheaves of $\ZZ(n)$ on $\Sm_\QQ$ in degrees $i \geq n$ are, \cite[Thm.10.1]{Blo86}, \cite{NS90}, \cite[\S3]{Tot92}, \cite[Thm.1.1]{Ker09},
\[
\sH^i(\ZZ(n)) = \left \{ \begin{array}{ll}
\sK_n^M & i = n \\
0 & i > n
\end{array} \right . 
\]
we have canonical morphisms 
\begin{equation} \label{equa:ZZnsKn}
\ZZ(n) \to \sK_n^M[-n]
\end{equation}
in $D(\PSh(\Sm_\QQ))$. On the other hand, for any $a_1, \dots, a_n \in A^*$, the form $\tfrac{da_1}{a_1}  \wedge \dots \wedge \tfrac{da_n}{a_n}$ lies in the kernel of the differential $\Omega^n_A \to \Omega^{n+1}_A$. So the morphism $\sK_n^M \to \Omega^n$ induces a morphism 
\begin{equation} \label{equa:SKOmega}
\sK_n^M[-n] \to \Omega^{\geq n}.	
\end{equation}
Composing Eq.\eqref{equa:ZZnsKn} with Eq.\eqref{equa:SKOmega} gives the desired morphism Eq.\eqref{equa:ZZnOmega}.
\end{cons}

\begin{cons}
We could use a similar procedure as Construction~\ref{cons:ZZnOmega} to construct a morphism Eq.\eqref{equa:ZZnnu} in $D(\PSh(\Sm_{\FF_p}))$, but we use instead a slightly different method. The system $(ν_\bullet^n)$ is a system of presheaves on $\SmCor_k$, \cite[Thm.8.3]{GL00}. The canonical element $\tfrac{dt_1}{t_1} \wedge \dots \wedge \tfrac{dt_n}{t_n}$ in $\lim_r ν_r^n(\GG_m^{\times n})$, where the $i$th $\GG_m$ is $\Spec(k[t_i, t_i^{-1}])$, defines a morphism 
\[ \hom_{\SmCor_k}(-, \GG_m^{\times n}) \to \varprojlim_r ν_r^n \]
in $\PSh(\SmCor_k)$. 
Now each $ν_r^n$ is $\AA^1$-invariant, \cite[Thm.8.3]{GL00}, \cite[Thm.4.27]{VoeCohThe}, and a Zariski sheaf, so we obtain a unique factorisation\footnote{Here we implicitly use the fact that $\lim_r ν_r^n = R\lim_r ν_r^n$.}
\begin{equation} \label{equa:LLGMnu}
 L^{\Zar}L^{\AA^1}\hom_{\SmCor_k}(-, \GG_m^{\times n}) \to \varprojlim_r ν_r^n.
\end{equation}
Finally, $\ZZ(n)[n]$ is a direct summand of $L^{\Zar}L^{\AA^1}\hom_{\SmCor_k}(-, \GG_m^{\times n})$.%
\footnote{By the Mayer-Vietoris triangle for the standard covering of $\PP^1$, we get $\Cone(M(\AA^1) \to M(\PP^1)) \cong \Cone(M(\GG_m) \to M(\AA^1))$. Since $\AA^1$ is contractible, this says $\ZZ(1)[2] = \tilde{M}(\PP^1) \cong \tilde{M}(\GG_m)[1]$ where $\tilde{M}(X)$ means the reduced motive of a scheme with respect to a rational point $x \in X$, i.e., $M(X) \cong \tilde{M}(X) \oplus M(x)$. Using the fact that $\DM_k^{\eff}$ has a unique tensor structure satisfying $M(X) \otimes M(Y) \cong M(X \times Y)$, we get isomorphisms $\ZZ(n)[2n] \cong \tilde{M}(\GG_m)[1]^{\otimes n} \subseteq M(\GG_m^{\times n})[n]$.}
So Eq.\eqref{equa:LLGMnu} induces a morphism Eq.\eqref{equa:ZZnnu}.
\end{cons}

\begin{defi}[Elmanto, Morrow, {\cite{EM23}}] \label{defi:EM}
On $\Sch_k$ define
\begin{equation*}
\ZZ(n)^{\EM}
= 
\Cone \biggl (\ZZ(n)^{\cdh} \oplus \ZZ(n)^{TC} \to L^{\cdh}(\ZZ(n)^{TC}) \biggr )[-1]	
\end{equation*}
where
\[ \ZZ(n)^{\cdh} := L^{\cdh}L^{\Kan}\ZZ(n), \]
cf.Def.\ref{defi:Znprocdh}, and $\ZZ(n)^{TC}$ is the 
Zariski local object whose values on affine $k$-schemes $\Spec(A)$ is
\begin{equation*}
\ZZ(n)^{TC}(A) = 
\left \{ 
\begin{array}{cc}
R\lim_m L^{\Kan} ((\Omega_A^{\geq n})^{\leq m}), &k = \QQ, \\
R\lim_r L^{\Kan} R\Gamma_{\et}(A, ν_r^n)[-n], &k = \FF_p.
\end{array}
\right .	
\end{equation*}
The morphisms
\begin{equation*}
\ZZ(n)^{\cdh} \to L^{\cdh}(\ZZ(n)^{TC})
\end{equation*}
are the ones induced by Eq.\eqref{equa:ZZnOmega} and Eq.\eqref{equa:ZZnnu} respectively.
\end{defi}

In \cite[Thm.8.2]{EM23} it is established that $\ZZ(n)^{\EM}$ satisfies procdh excision, and therefore by Proposition~\ref{prop:MVprocdh} has procdh descent. Consequently, one obtains canonical comparison maps
\[ \ZZ(n)^{\procdh} \to \ZZ(n)^{\EM}. \]

The following comparison theorem was obtained through joint discussion with Morrow. The proof relies heavily on the main results of \cite{EM23}.

\begin{theo}[K., Saito, {\cite[Cor.1.11]{KS24}}, Elmanto, Morrow, \cite{EM23}] \label{theo:procdhEM}
For $k = \FF_p$ or $\QQ$ and $Y \in \Sch_k$ there are equivalences 
\begin{equation} \label{equa:comparison}
\ZZ(n)^{\procdh}(Y) \cong \ZZ(n)^{\EM}(Y).
\end{equation}
\end{theo}

The strategy for the comparison in Theorem~\ref{theo:procdhEM} is straight-forward. Both $\ZZ(n)^{\procdh}$ and $\ZZ(n)^{\EM}$ are procdh sheaves, the former by definition and the latter by \cite[Thm.8.2]{EM23} and Proposition~\ref{prop:MVprocdh}.  So if the procdh topos has enough points, it is enough to compare them on procdh local rings.

Identification of local rings for the procdh topology is something we did not have room for here. Local rings for the cdh topology are henselian valuation rings, \cite{GL01}, \cite{GK15}. That is, henselian valuation rings are to the cdh topology as strictly henselian local rings are to the étale topology. Local rings for the procdh topology are ``nice'' nilpotent thickenings of henselian valuation rings. Explicitly, a ring is \emph{procdh local} if it is of the form $\OO \times_K A$ for some local ring $A$ of Krull dimension zero, and valuation ring $\OO \subseteq K = A / \m$ of its residue field, \cite[Prop.3.5]{KS24}. For example, $\{a + b ε \in \QQ_p[ε] / ε^2\ |\ a \in \ZZ_p\}$ is a procdh local ring. Since procdh coverings are not finite, it doesn't follow from Deligne's completeness theorem that the procdh topos has enough points. None-the-less, the procdh topos \emph{does} have enough points, \cite[Thm.1.5]{KS24}. Describing the values of $\ZZ(n)^{\EM}$ on procdh local rings is done in \cite{EM23}.


\begin{rema}\ \label{rema:ZEMrem}
\begin{enumerate}
 \item 	The complexes $\ZZ(n)^{\procdh}$ and $\ZZ(n)^{\EM}$ are actually defined for qcqs schemes over a field. Theorem~\ref{theo:procdhEM} holds for all Noetherian schemes over $\FF_q$ or $\QQ$. 

 \item The construction of the motivic complex and development of the theory in mixed characteristic has been worked out now by Tess Bouis. In particular, Bouis has proved Theorem~\ref{theo:procdhEM} for general Noetherian schemes (also in mixed characteristic), \cite{Bou24}.

 \item  The Atiyah-Hirzebruch spectral sequence with $\ZZ(n)^{\EM}$ is proved in \cite{EM23} for qcqs schemes over a field, completely independently from \cite{KS24}. For $\ZZ(n)^{\procdh}$ both the complex and the spectral sequence can be defined on qcqs schemes over $k$, but convergence to algebraic $K$-theory is only known for Noetherian schemes. In general, the target is%
\footnote{Here we abuse notation and use $L^{\procdh}$ also for the sheafification $\PSh(\Sch^{\qcqs}_k, \Spt) \to \Shv_{\procdh}(\Sch^{\qcqs}_k, \Spt)$ of presheaves of spectra on the category $\Sch^{\qcqs}_k$ of all qcqs $k$-schemes.} %
$L^{\procdh}K$ and convergence is unclear.

 \item As a guiding problem we proposed the question of extending Voevodsky/Bloch's motivic complexes $\ZZ(n)$ beyond smooth varieties. One of the central theorems in \cite{EM23} is that the restriction of $\ZZ(n)^{\EM}$ to smooth schemes recovers $\ZZ(n)$. 
\end{enumerate}
\end{rema}

\begin{rema} \label{rema:cdhZ}
The relationship of $\ZZ(n)^{\cdh} = L^{\cdh}L^{\Kan}\ZZ(n)$ to the category $DM_{k,\cdh}^{\eff}$ is a little subtle without some kind of resolution of singularities. For $X \in \Sm_k$, the identification $\hom_{\SmCor_k}(-, X) \stackrel{=}{\to} \hom_{\Cor_k}(-, X)|_{\SmCor_k}$ induces a comparison 
\[ L^{\cdh}L^{\Kan}L^{\Zar}L^{\AA^1}\hom_{\SmCor_k}(-, X) \to L^{\cdh}L^{\AA^1}\hom_{\Cor_k}(-, X). \]
With some kind of resolution of singularities and appropriate coefficients, this is an equivalence.%
\footnote{For example, with $\ZZ$-coefficients, assuming $k$ satisfies resolution of singularities this is in \cite[\S 4]{Voe00}. With $\QQ$-coefficients it is observed in \cite[Intro.]{SV00} that de Jong's theorem on alterations can be used as a replacement for resolution of singularities. With $\ZZ[1/p]$-coefficients it is proved in \cite{Kel17} that Gabber's theorem on alterations can be used as a replacement for resolution of singularities.} %
However, without assuming resolution of singularities, comparing the $p$-torsion is still open. Consequently, $\ZZ(n)^{\cdh}$ as defined above and the complex $\ZZ(n)^{DM_{\cdh}}$ represented in $\DM_{k,\cdh}$ might have different $p$-completions. (This would disprove resolution of singularities in positive characteristic.) In forthcoming work of Bachmann, Elmanto, and Morrow, it is shown that:
 The $\ZZ(n)^{\cdh}$ are the graded pieces of the slice filtration on $KH(X)$, when $X$ is any qcqs $k$-scheme.

\end{rema}

\section{But what is a motive, though?} \label{sec:5}

We began the talk with a category of motives, and then ended up with a motivic cohomology. We observed in Remark~\ref{rema:cdhA1} that correspondences as they are defined in \cite{SVRelCyc}, the first step to building $\DM_k^{\eff}$, can't see nilpotents. There is however a category, not of motives, but of stable motivic homotopy types which can see nilpotents, namely the category of motivic spectra $\MS_S$ of Annala and Iwasa, \cite{AI22}, further developed by them with Hoyois, \cite{AHI24}.

%
%
%
It is currently an open question whether $S \mapsto \MS_S$ satisfies procdh descent as a presheaf of categories. However, a projective bundle formula and elementary blow-up excision for $\ZZ(j)^{\EM}$ are established in \cite{EM23}, which means that $\ZZ(j)^{\EM}$ is represented in $\MS_S$, for $S$ any qcqs $k$-scheme. One could try considering a category of modules over it in analogy to \cite{Rob87}, \cite{RO06} and \cite{HKO17}.

\bibliographystyle{alpha}

\newcommand{\etalchar}[1]{$^{#1}$}

\end{document}